\documentclass{lmcs} 
\pdfoutput=1
\usepackage[utf8]{inputenc}

\usepackage{lastpage}
\lmcsdoi{22}{4}{11}
\lmcsheading{}{\pageref{LastPage}}{}{}%
{Jul.~10,~2025}{Oct.~27,~2026}{}


\usepackage{hyperref}
\usepackage{amssymb}
\usepackage{tikz}
\usetikzlibrary{cd}
\usepackage{graphicx}
\usepackage{mathrsfs}
\usepackage{comment}
\theoremstyle{plain} 

\newcommand{\rom}[1]{\rm{\uppercase\expandafter{\romannumeral #1}}}

\begin{document}
		
		\title[Commutation of Power Constructions in Well-filtered Dcpos]{Commutation of Smyth and Hoare Power Constructions in Well-filtered Dcpos}
		\thanks{This work is supported by the National Natural Science Foundation of China (No.~12231007, No.~5kj20250124 and No.~12601896) and the Natural Science Foundation of Shandong Province, China (No.~ZR2024QA195).}
		\thanks{Corresponding author: Qingguo Li, \texttt{liqingguoli@aliyun.com}.}
		
		\author[H.~Hou]{Huijun Hou\lmcsorcid{0000-0001-8474-0433}}[a]
		\author[Q.~Li]{Qingguo Li\lmcsorcid{0000-0002-4314-7489}}[b]
		
		\address{School of Mathematical Sciences, Qufu Normal University, Qufu, Shandong, 273165, China}
		\email{houhuijun2021@163.com}
		
		\address{School of Mathematics, Hunan University, Changsha, Hunan, 410082, China}
		\email{liqingguoli@aliyun.com}
		
		
		\begin{abstract}
			Prior work \cite{heck1991} established a commutativity result for the Hoare power construction and a modified version of the Smyth power construction consisting of strongly compact sets, which is defined for $\mathcal U_{S}$-admitting dcpos, where $\mathcal U_{S}$-admissability is well-filteredness with compact sets replaced by strongly compact sets.
			In this paper, we consider the  Hoare power construction  $\mathcal H$ and the Smyth power construction  $\mathcal Q$  on the category $\mathbf{WF}$ of  well-filtered dcpos with Scott-continuous maps.
			Actually, the functors $\mathcal H$ and  $\mathcal Q$ can be extended to monads.
			We prove that $\mathcal H$ and $\mathcal Q$ commute, that is, $\mathcal {HQ}(L)$ is isomorphic to $\mathcal {QH}(L)$ for a well-filtered dcpo $L$,  if and only if $L$ satisfies a property similar to consonance that we call $(\mathrm{KC})$ and the Scott topology coincides with the upper Vietoris topology on $\mathcal Q(L)$. We also investigate the Eilenberg-Moore category of the monad composed by $\mathcal H$ and $\mathcal Q$ under a distributive law $\phi$ on $\mathbf{WF}$ and characterize it  to be a  subcategory of the category $\mathbf{Frm}$, which is composed of all frames and all frame homomorphisms.
		\end{abstract}
		
		\maketitle
		
\section{Introduction}

Non-deterministic choices are important semantic concepts  that offer new insights in designing more powerful programming languages. They were originally modeled by the power constructions on the category of domains. The most common of these are  the Hoare and the Smyth  powerdomains \cite{smyth}, which  can respectively serve as the denotational semantics for  angelic and demonic non-deterministic choices. Notably, the Hoare and the Smyth  powerdomains also give rise to monads. In 1991, E. Moggi showed that, within the category of sets, the collection of all finite sets could serve as denotations for nondeterministic programming languages  \cite{moggi}. He further proved that this structure forms a monad on the category of sets, thereby linking power constructions with monads. Over time, monads have evolved into a fundamental tool for constructing denotational semantic models in functional programming languages.

In 1969, J. Beck proposed the notion of the distributive law between two monads defined on the same category \cite{beck}. 
Such a distributive law   ensures that the composition of these two monads remains a monad,  providing a new approach to establishing denotational semantic models for functional programming languages. Moreover, the two composite monads obtained under different composition orders may differ; if they coincide, i.e., the monad functors induce isomorphic mappings on objects, then the monads are said to obey a commutative law. Such commutativity ensures that the order of program execution  becomes immaterial in computation. This practical significance raises two fundamental problems: (1) Across which categories does a distributive law exist between  the two denotational semantic models for nondeterministic programming languages -- the Hoare power monad and the  Smyth power monad? (2) Within a specified category, under what conditions do these two monads commute?

In domain theory, numerous scholars have conducted in-depth investigations into the above two problems. In 1990, K. Flannery and J. Martin employed information systems to demonstrate the commutativity between the Hoare power monad and the Smyth power monad in the category of bounded complete algebraic domains \cite{fm1990}. Subsequently, R. Heckmann used an entirely distinct method to prove that for  the category of dcpos, the commutativity of Hoare and Smyth power monads is guaranteed for variants of the two monads that are defined as the free inflationary semilattice and the free deflationary semilattice, respectively \cite{heckdcpo}. Within the category of domains, the Hoare powerdomain induces the free up-complete $\vee$-semilattices, while the Smyth powerdomain yields the free meet-continuous semilattices \cite{clad}. Here, the up-complete $\vee$-semilattices are actually inflationary semilattices, and the  meet-continuous semilattices are deflationary semilattices. Thus the Hoare power monad and the Smyth power monad commute in the category of domains.

Note that  many categories   have the category of domains as a full subcategory. The versions of Hoare and Smyth power monads considered in \cite{heckdcpo}, i.e. defined as free constructions, are well-defined for all dcpos and do commute
on the category of all dcpos. 
As shown by A. Schalk in \cite{as}, we know  the Hoare power construction defined by the Scott closed subsets does yield the free
unital inflationary semilattice in $\mathbf{DCPO}$;  however,  it is not  clear what constructions can  give a free deflationary semilattice in $\mathbf{DCPO}$.  	We know that the well-filteredness guarantees that the collection of all compact saturated subsets with the reverse inclusion order is a dcpo \cite{clad}.
So the issue of  restricting 
to well-filtered dcpos (called $\mathcal U_K$-admitting dcpos in \cite{heck1991})  in order to obtain a
well-behaved Smyth construction  has attracted increasing attention. Concretely, R. Heckmann in \cite{heck1991} investigated the commutativity between the Hoare power construction $\mathcal L$  and a modified Smyth power construction $\mathcal U_S$ on the category of $\mathcal U_S$-admitting dcpos \cite{heck1991}. In this work,  the modified Smyth power construction is determined by all strongly compact  subsets, which  differs from  the traditional Smyth power construction based on compact subsets. Thus, the   $\mathcal  U_{S}$-admissibility is strictly weaker than  well-filteredness. R. Heckmann introduced the property of being  $\mathcal U_S$-conformal and showed that for a dcpo $L$, $L$ is $\mathcal U_S$-conformal if and only if $\mathcal L(\mathcal U_S(L))$ is isomorphic to $\mathcal U_S(\mathcal L(L))$, i.e.,  established the sufficient and  necessary conditions  for the commutativity. Meanwhile, they also pointed out that they could not show
an analogous commutativity result for $\mathcal L$ and the Smyth power construction based on compact subsets.

In this paper, we focus on the category of well-filtered dcpos and  investigate the commutativity  between the associated Hoare and Smyth power monads. In Section 3, we  define the Hoare and Smyth constructions  for a well-filtered dcpo $L$  as follows: The Hoare construction $\mathcal H$ is the set of all Scott closed subsets of $L$, ordered by inclusion; the Smyth construction $\mathcal Q$ is the set  of all  Scott compact saturated subsets of $L$, ordered by  reverse inclusion order. We explicitly show that $\mathcal Q$ does not yield a free construction for the category  $\mathbf {WF}$ of  well-filtered dcpos with Scott-continuous maps, which confirms the significance and necessity of our research.  In Section 4, we  introduce a  property called $(\mathrm{KC})$, which is 
 weaker than the  consonance  in terms of  posets. Based on this, we proved that $\mathcal H$ and $\mathcal Q$ commute in the sense of that $\mathcal{HQ}(L)$ is isomorphic to $\mathcal{QH}(L)$ under a certain function, if and only if $L$ satisfies the property $(\mathrm{KC})$ and the Scott topology coincides with the upper Vietoris topology on $\mathcal Q(L)$. It is necessary to point out that in \cite{brecht},  M. de Brecht and  T. Kawai also provided a Smyth construction via compact saturated subsets and investigated the commutativity between the Smyth construction and the Hoare construction. However, our work differs from theirs,   since they
 consider the category of  topological spaces while we
 explore that of well-filtered dcpos.
 
 Notably, unlike R. Heckmann's modified Smyth power construction $\mathcal U_{S}$, $\mathcal Q$  induces  the Smyth power monad on  $\mathbf{WF}$. This relies crucially on  Xu's result  \cite{Xu} that the Smyth power space endowed with the Scott topology of any well-filtered space remains well-filtered. Meanwhile, since  every  complete lattice  with the Scott topology is well-filtered \cite[Corollary 3.2]{Xi}, we conclude  that $\mathcal H$  also gives rise to a monad on  $\mathbf{WF}$. Consequently,  in Section 5, we show that there is a distributive law between the two monads $\mathcal Q$ and $\mathcal H$,  which ensures that the composition  $\mathcal{QH}$  remains  a monad on $\mathbf{WF}$. Finally, motivated by the growing interest in Eilenberg-Moore algebras, we characterize the Eilenberg-Moore category of  $\mathcal{QH}$  as a subcategory of \textbf{Frm} -- the category of frames and frame homomorphisms.

\section{Preliminaries}

In this section, we briefly introduce some basic concepts and notations that will be used in this paper.
One can  refer to \cite{ct,clad, nht2} for more details.

Given a poset $P$, $A\subseteq P$ is an \emph{upper set} if $A = {\uparrow} A$, where ${\uparrow} A = \{x\in P: x\geq a\  \mathrm{for\  some}\ a\in A\}$. The  \emph{lower set} is defined dually. A subset $D$ of $P$ is \emph {directed}  if it is nonempty and every finite subset  of $D$ has an upper bound in $D$. If each directed subset in  $P$ has a supremum, we call $P$  a {\em directed
complete partially ordered set (dcpo)}.  A subset $U$ of $P$ is called \emph{Scott open} iff it is an upper set and for each directed subset $D$ whose supremum exists, written as $\sup D$, $\sup D\in U$ implies $D\cap U\neq \emptyset$. The complement of a Scott open set is called \emph{Scott closed}. The collection of all Scott open subsets, denoted by $\sigma(P)$, forms a topology  called \emph{Scott topology} and in general, we write $\Sigma P = (P, \sigma(P))$. Besides, we use $\Gamma(P)$ to denote the set of all Scott closed subsets of $P$.  A continuous map between two posets equipped with Scott topology is called \emph{Scott-continuous}. Generally, $\mathbf{DCPO}$ is used to denote the category of all dcpos with Scott-continuous mappings.

Let $X$ be a topological space. We write $\mathcal O(X)$ for the collection  of  open subsets of $X$.  When  $X$ is a $T_{0}$ space, the  \emph{specialization preorder} $\leq$ defined by $x\leq y$ iff $x\in cl(\{y\})$ is a partial order and   is called \emph{specialization order}.  A subset of $X$ is \emph{saturated} if it is the intersection of open sets, or equivalently, it is an upper set in the order of specialization.  $X$ is called \emph{well-filtered} if for each filter basis $\mathcal C$ of compact saturated subsets and any open subset $U\in \mathcal O(X)$ with $\bigcap \mathcal C\subseteq U$, $K\subseteq U$ for some $K\in \mathcal C$. We usually say a dcpo $P$ is well-filtered if it, equipped with the Scott topology, is a well-filtered space. Furthermore, a space $X$ is called \emph{coherent} if the intersection of any two compact saturated subsets of $X$ is also compact.


In this paper, we
denote the set of all  compact saturated subsets (including the empty set) and that of all closed subsets of a space $X$ by $\mathcal Q(X)$ and $\Gamma(X)$, respectively. The topology on $\mathcal Q(X)$ generated by $\{\square U: U\in \mathcal O(X)\}$ as a base, where  $\Box U = \{K\in  \mathcal Q(X): K\subseteq U\}$, is called the \emph{upper Vietoris topology} and denoted by $\upsilon(\mathcal Q(X))$. The upper topology on $\Gamma(X)$ coincides with the \emph{lower Vietoris topology}, which is written as $\nu(\Gamma(X))$ and  generated by $\{\lozenge U: U\in \mathcal O(X)\}$ as a subbase, where $\lozenge U = \{A\in \Gamma(X): A\cap U\neq\emptyset\}$.

	Regarding the  inflationary (deflationary) semilattices, there is a formal definition  given in \cite{clad}, and a concrete description  provided by A. Schalk in \cite{as}. Combining these two forms, we give the following definition:

\begin{defi}
An inflationary (deflationary) semilattice is a dcpo equipped with a Scott-continuous binary operation
$\uplus$ that is commutative, associative, idempotent, and satisfies $x \leq x \uplus y \;(x \geq x \uplus y)$.
In addition, an inflationary (deflationary) semilattice is unital if $\uplus$ has a neutral element.
\end{defi}

\begin{defiC}[\cite{ct}]
	A \emph{monad} $\mathcal T$ in a category $\mathbf C$ is a triple consisting of an endofunctor $\mathcal T: \mathbf C\rightarrow\mathbf C$ and two natural transformations: the \emph{unit}	$\eta: Id_{\mathbf C}\rightarrow \mathcal T$  and the \emph{multiplication} $\mu: \mathcal T^{2}\rightarrow \mathcal T$ 	such
that the following diagrams commute.
	\begin{figure}[ht]
	\centering
	\begin{tikzpicture}[line width=0.6pt,scale=0.7]
		\node[right] at(-5,1.1)(a)  { $\mathcal T$};
		\node[right]  at(-2,1.1) (b)  { $\mathcal T^{2}$};
		\node[right] at(1,1.1) (c){ $\mathcal T$};
		\node[right]  at(-1.9,-1.1) (d)  { $\mathcal T$};
		\draw[->](a)--(b);
		\draw[->](c)--(b);
		\draw[->](a)--(d);
		\draw[->](c)--(d);
			\draw[->](b)--(d);
		\node[right]  at(-3.5,1.3)  { $\eta_{\mathcal T}$};
		\node[right] at(-0.5,1.3)  { $\mathcal T\eta$};
		\node[right] at(-1.6,0)  {$\mu$};
		\node[right]  at(-3.9,-0.2)  {=};
		\node[right]  at(-0.1,-0.2)  { =};
		
		\node[right] at(5,1.1)(d)  {$\mathcal T^{3}$};
		\node[right] at (8.5,1.1) (e){ $\mathcal T^{2}$};	
		\node[right]  at(5.05,-1.1) (g)  { $\mathcal T^{2}$};
		\node[right]  at(8.65,-1.1) (h)  { $\mathcal T$};
		\draw[->](d)--(e);
		\draw[->](e)--(h);
		\draw[->](d)--(g);
		\draw[->](g)--(h);
		\node[right] at(6.6,1.4)  { $\mathcal T\mu$};
		\node[right] at(4.5,0)  { $\mu_ {\mathcal T}$};
		\node[right] at(8.9,0)  { $\mu$};
		\node[right] at(7,-1.4)  { $\mu$};
		
	\end{tikzpicture}
\end{figure}
\end{defiC}

\begin{defiC}[\cite{ct}]
	Let $\mathcal T = (\mathcal T, \eta, \mu)$ be a monad on the category $\mathbf{C}$.
	\begin{enumerate}[(1)]
		\item A  pair $(C, \alpha)$ is a  $\mathcal{T}$$\emph{-algebra}$ (or \emph{Eilenberg-Moore algebra} of  $\mathcal{T}$)  if $\alpha: \mathcal{T}C\rightarrow C$ is a morphism in $\mathbf{C}$ that satisfies:
		
		\begin{center}
			($\mathbf{Associativity}$) $\alpha\circ\mu_{C} = \alpha\circ \mathcal{T}\alpha$,\qquad  ($\mathbf{Unit\; law}$) $\alpha\circ \eta_{C} = id_{C}$.
		\end{center}
		In this case, $\alpha$ is called a \emph{structure map}.
		\item Let $(A, \alpha_{A})$ and $(B, \alpha_{B})$ be $\mathcal{T}$-algebras. A \emph{$\mathcal T$-algebra homomorphism} from $(A, \alpha_{A})$ to $(B, \alpha_{B})$ is an arrow $h: A\rightarrow B$ satisfying:
		\begin{center}
			$(\mathbf{Homomorphism\;law})\ h\circ \alpha_{A} = \alpha_{B}\circ \mathcal{T}h$.
		\end{center}
	\end{enumerate}

\noindent
The category composed  of  $\mathcal{T}$-algebras and $\mathcal{T}$-algebra homomorphisms is  called \emph{Eilenberg-Moore category} of $\mathcal{T}$ over
$\mathbf C$.
\end{defiC}

\begin{prop} \label{mp}
	Let $\mathcal T = (\mathcal T, \eta, \mu)$ be a monad on the category $\mathbf{C}$. Then the following properties hold:
	\begin{enumerate}[(1)]
	\item  	For every object $C$ in $\mathbf{C}$, $(\mathcal TC, \mu_C)$ is a $\mathcal T$-algebra.
	\item 	 If $(C, \alpha)$ is a $\mathcal T$-algebra, then $\alpha: (\mathcal TC, \mu_C) \rightarrow (C, \alpha)$ is a $\mathcal T$-algebra homomorphism.
	\item  If $h_1, h_2:  (\mathcal TC, \mu_C)\rightarrow  (C', \alpha)$ are $\mathcal T$-algebra homomorphisms with $h_1 \circ \eta_C = h_2 \circ \eta_C$, then $h_1 = h_2$.

	\end{enumerate}
		\begin{proof}
		(1) and (2) are clear from the definitions, more precisely, one can refer to \cite{maclane}. Part (3) holds because $h_i = h_i \circ \mu_C \circ \mathcal T\eta_C = \alpha \circ \mathcal Th_i \circ \mathcal T\eta_C = \alpha \circ \mathcal T (h_i \circ \eta_C) $.
		\qedhere
	\end{proof}
\end{prop}

\section{The Hoare and Smyth power constructions}
\subsection{The Hoare power construction}
\begin{defi}
The \emph{Hoare power construction} $\mathcal H(L)$ for a well-filtered dcpo $L$ is  the set of all Scott closed subsets of $L$ ordered by inclusion.
\end{defi}

\noindent
We now focus on constructing mutually adjoint functors, through which the composite functor  induces the monad $\mathcal{H}$ on $\mathbf{WF}$. Clearly, for each well-filtered dcpo $L$, its Hoare power construction is a complete lattice.
 Let $\mathbf{CL}$  denote the category of complete lattices with morphisms preserving arbitrary suprema.
 
Consider the morphisms in $\mathbf{WF}$. For any Scott-continuous mapping $f: L \rightarrow M$  between well-filtered dcpos $L$ and $M$, we define $\mathcal H(f): \mathcal H(L)\rightarrow \mathcal H(M)$ as
\begin{center}
 $\mathcal H(f)(A) = cl(f(A))$,	for any $A\in \mathcal H(L)$.
\end{center}
Let $\{A_{i}: i\in I\}$ be an arbitrary subset of $\mathcal H(L)$. Since $\mathcal H(f)$   is clearly order-preserving,  $\sup_{i\in I}\mathcal H(f)(A_{i})\subseteq \mathcal H(f)(\sup_{i\in I}A_{i})$ holds. Meanwhile, 
\begin{center}
	$\mathcal H(f)(\sup_{i\in I}A_{i}) = \mathcal H(f)(cl(\bigcup_{i\in I}A_{i})) \subseteq cl(f(\bigcup_{i\in I}A_{i})) = cl(\bigcup_{i\in I}f(A_{i})) = \sup_{i\in I}\mathcal H(f)(A_{i})$,
\end{center}
where the  symbol $``\subseteq"$ is guaranteed by the continuity of $f$.
Thus we conclude that $\mathcal H(f)$ preserves arbitrary suprema.
Additionally, one can verify that $\mathcal H$ preserves the identities and the compositions of any two morphisms. Along with the above, we can obtain the following result.

\begin{lem}
	$\mathcal H$ is a functor from $\mathbf{WF}$ to $\mathbf{CL}$.
\end{lem}

\noindent
It follows from \cite[Corollary 3.2]{Xi} that each complete lattice endowed with the Scott topology is well-filtered. So   there is a forgetful functor  from $\mathbf{CL}$ to $\mathbf{WF}$.  A. Schalk has proved that the functor $P_D^b: \mathbf{DCPO}\rightarrow \mathbf{CL}$ is left adjoint to the forgetful functor $U: \mathbf{CL} \rightarrow\mathbf{DCPO}$  (see \cite[Theorem 6.2]{as}). Note that if the functor $P_D^b$ is restricted from the category $\mathbf{DCPO}$ to its full subcategory $\mathbf{WF}$, it is just the functor $\mathcal H$. This immediately  yields  the following result. 

\begin{prop}\label{hl}
	The functor $\mathcal H: \mathbf{WF}\rightarrow \mathbf{CL}$ is left adjoint to the forgetful functor $U: \mathbf{CL}\rightarrow\mathbf{WF}$. 
\end{prop}

\noindent
Since each pair of adjunctions determines a monad,  \cite[Proposition 4.2.1]{ct},	we deduce that the triple $(U\circ \mathcal H, \eta, U\varepsilon_{\mathcal H})$,  where $\eta$ and $\varepsilon$ are the unit and counit respectively, turns into a monad on $\mathbf{WF}$. For ease of expressions, we will abbreviate $U\circ\mathcal H$ as $\mathcal H$. Meanwhile,  when $U\varepsilon_{\mathcal H}$ acts on a well-filtered $\mathrm{dcpo}$ $L$, one can calculate $U\varepsilon_{\mathcal H}(L)$ as the natural transformation $\mu_L$ from  $\mathcal{H}(\mathcal{H}(L))$ to $\mathcal{H}(L)$  that maps $\mathcal{A}$ to $\sup_{\mathcal{H}(L)}\mathcal A$.

\begin{lem}
	Let $L$ be a dcpo and $\mathcal A$ a Scott closed subset of $\mathcal H(L)$. Then $\bigcup\mathcal A$ is Scott closed.
\end{lem}
\begin{proof}
It is similar to that of   \cite[Proposition 2.2]{Ho}.
\qedhere
\end{proof}
\noindent
It follows from this lemma  that $\sup_{\mathcal{H}(L)}\mathcal A = \bigcup\mathcal A$. So we  conclude that
\begin{thm}
	The endofunctor $\mathcal{H}$ together with the unit $\eta: Id\rightarrow \mathcal H$ and the multiplication $\mu: \mathcal H^{2}\rightarrow \mathcal H$ forms a monad on $\mathbf{WF}$, where $\eta$ and $\mu$ are defined concretely for a well-filtered dcpo $L$ as follows:
	\begin{center}
		 $\eta_{L}(x) = {\downarrow }x$, for any $ x\in L$,
	\end{center}
	and
	\begin{center}
 $\mu_{L}(\mathcal{A}) = \bigcup \mathcal{A}$, 	for any	$\mathcal{A}\in \mathcal{H}(\mathcal{H}(L))$.
	\end{center}
\end{thm}

\noindent
Since $\mathbf{CL}$ is the category of unital inflationary semilattices for $\mathbf{WF}$, and Proposition \ref{hl} shows that  the Hoare power construction $\mathcal H(L)$ has the universal property for every well-filtered dcpo $L$,  we  obtain the following result.
\begin{cor}
	$\mathcal H$ gives the free unital inflationary semilattice for $\mathbf{WF}$.
\end{cor}

\begin{lem}\label{ha}
	Let $L, M$ be  well-filtered dcpos. Then
\begin{enumerate}[(i)]
	\item $(L, \alpha)$ is an $\mathcal H$-algebra if and only if $L$ is a complete lattice.
	\item $f: (L, \alpha)\rightarrow (M, \beta)$ is an $\mathcal H$-algebra homomorphism if and only if $f$ is a complete lattice homomorphism.
\end{enumerate}
\begin{proof}
Based on their definitions, it is easy to verify the claims.	
\qedhere
\end{proof}
\end{lem}

\begin{cor}\label{halgebra}
	The Eilenberg-Moore category of $\mathcal H$ over $\mathbf{WF}$ is just $\mathbf{CL}$.
\end{cor}

\subsection{The Smyth power construction}

Xu et al. proved that for a  well-filtered space $X$, $\Sigma \mathcal Q(X)$ remains well-filtered \cite[Theorem 6.5]{Xu}, which indicates that for each
 well-filtered dcpo $L$, $\mathcal Q(L)$ with the Scott topology is  well-filtered, that is,  $\mathcal Q(L)$ is an object in $\mathbf{WF}$. This fact  gives us a well-defined Smyth power construction.

\begin{defi}
	The \emph{Smyth power construction}  for a well-filtered dcpo $L$ is the set $\mathcal Q(L)$ of all  Scott compact saturated subsets  ordered by  reverse inclusion.
\end{defi}

\noindent
Recall that, inspired by the facts that for a sober space $X$, the set $\mathcal Q(X)$ ordered by reverse inclusion is isomorphic to the set of Scott open filters in $\mathcal O(X)$ ordered by inclusion, and for a general dcpo $L$, the set $\mathcal Q(L)$ of all compact saturated subsets  ordered by  reverse inclusion is no longer a dcpo, A. Schalk focused on the set of all Scott open filters $OFilt(\sigma(L))$ and proved that $OFilt(\sigma(L))$ ordered by inclusion remains a dcpo. Concerning this, she 
	gave an endofunctor on the category $\mathbf{DCPO}$  and obtained a monad (one can refer to \cite[Lemma 7.5, Theorem 7.7]{as}). 
	Fortunately, for a well-filtered dcpo $L$, the set of all Scott compact saturated subsets ordered by reverse inclusion remains a well-filtered dcpo as shown in \cite{Xu}. Next we will construct an endofunctor $\mathcal Q$  on $\mathbf{WF}$ and verify that it is also a  monad.

 \begin{lemC}[{\cite[Theorem 3.9]{heck1991}}]
 	Let $f: L\rightarrow M$ be a Scott-continuous mapping between well-filtered dcpos $L$ and $M$. Then $\mathcal Q(f): \mathcal Q(L)\rightarrow \mathcal Q(M)$  defined as
 	\begin{center}
 	 $\mathcal Q(f)(K) = {\uparrow}f(K)$, 	for any $K\in \mathcal Q(L)$
 	\end{center}
 	is Scott-continuous.
 \end{lemC}

\noindent
One can easily verify that $\mathcal Q$ preserves the identities and the compositions of any two morphisms.  So $\mathcal Q$ is an endofunctor on $\mathbf{WF}$. 
\begin{lem}
Let $L$ be a well-filtered dcpo and $\mathcal K$  a Scott compact saturated  subset of $\mathcal Q(L)$. Then $\bigcup\mathcal K$ is compact saturated in $\Sigma L$.	
	\begin{proof}
		The saturation is clear, so we only need to consider the compactness. Let $\mathcal U$ be a directed subset of $\sigma(L)$ with $\bigcup\mathcal K\subseteq\bigcup\mathcal U$. Then for each $K\in \mathcal K$, its compactness yields the existence of $U_{K}\in \mathcal U$ such that $K\subseteq U_{K}$, that is, $K\in \Box U_{K}$.  Since the well-filteredness of $L$ implies that each $\Box U$ is Scott open in $\mathcal Q(L)$ and $\mathcal U$ is directed,   $\{\Box U: U\in \mathcal U\}$ forms a directed family of Scott open subsets of $\mathcal Q(L)$, and $\mathcal K\subseteq\bigcup\{\Box U: U\in \mathcal U\}$.  As $\mathcal K$ is compact, there exists a $U_{0}\in \mathcal U$ such that 
		 $\mathcal K\subseteq \Box U_{0}$, which immediately implies $\bigcup\mathcal K\subseteq U_{0}$. So $\bigcup\mathcal K$ is compact.
	\qedhere
	\end{proof}
\end{lem}

\begin{thm}
The endofunctor $\mathcal Q$ with the unit $\theta: Id\rightarrow \mathcal Q$ and the multiplication $\iota: \mathcal Q^2\rightarrow \mathcal Q$ forms a monad on $\mathbf{WF}$, where $\theta$ and $\iota$ are defined  as follows when they act on a well-filtered dcpo $L$:
\begin{center}
 $\theta_{L}(x) = {\uparrow}x$, 	for any $x\in L$,
\end{center}
and
\begin{center}
	 $\iota_{L}(\mathcal K) = \bigcup\mathcal K$, for any $\mathcal K\in \mathcal Q(\mathcal Q(L))$.
\end{center}
\begin{proof}
	It is  routine to check that $\mathcal Q$ is a monad.
	\qedhere
\end{proof}
\end{thm}
\noindent
It is noted that for a well-filtered dcpo $L$,  arbitrary finite infs and directed sups in $\mathcal Q(L)$ are actually the finite unions and filtered intersections. Thus $\mathcal Q(L)$ is a meet continuous semilattice, in other words, $\mathcal Q(L)$ is a deflationary semilattice in $\mathbf{WF}$.
\begin{lem}
Every $\mathcal Q$-algebra $M$ is a deflationary semilattice in $\mathbf{WF}$.
\begin{proof}
The definition of $\mathcal Q$-algebras reveals that $M$ is a retract of $\mathcal Q(M)$. Then by  \cite[Proposition 5.6]{as}, we know that $M$ is  a deflationary semilattice in $\mathbf{WF}$.
\qedhere
\end{proof}
\end{lem}

\noindent
Based on \cite[Lemma 4.4]{as} given by A. Schalk, we can derive the following results. 
\begin{lem}\label{inf}
	Let $(M, \alpha)$ be a $\mathcal Q$-algebra. Then  for each $K\in \mathcal Q(M)$, its infimum $\wedge K$ exists.
\end{lem}
\begin{lem}\label{qhom}
Each $\mathcal Q$-algebra homomorphism preserves infs of compact sets; hence, in particular finite infs.
\end{lem}
\begin{prop}\label{qalgebra}
The Eilenberg-Moore category of $\mathcal Q$ over $\mathbf{WF}$ is a subcategory of the category of deflationary semilattices in  $\mathbf{WF}$.
\end{prop}
\noindent
Unfortunately, we cannot accurately characterize what the  $\mathcal Q$-algebras are. Moreover, $\mathcal Q$ does not give the free deflationary semilattice for $\mathbf{WF}$ because there are deflationary
semilattices in $\mathbf{WF}$ that are not  $\mathcal Q$-algebras, as shown in the following example.

\begin{exa}\label{l}
	Let $\mathcal L$ be the dcpo given in Section 7.2.4 by A. Schalk in \cite{as}, that is,
	\begin{center}
		$\mathcal L = \{A\subseteq [0, 1]: A\neq\emptyset\ \mathrm{and} \ [0, 1]\setminus A \ \mathrm{is\ countable}\}$.
	\end{center}
It has been proved to be a deflationary semilattice in $\mathbf{DCPO}$. Besides, we can see that every nonempty finite subset of $\mathcal L$  has a sup, which is indeed the union of its elements. It then follows that
$\mathcal L$ is  an up-complete $\vee$-semilattice, this implies the well-filteredness of  $\Sigma \mathcal L$  witnessed by  \cite[Proposition 3.1]{Xi}. Thus $\mathcal L$ is a deflationary semilattice in $\mathbf{WF}$.  However, there exists a compact saturated subset
\begin{center}
	$K = {\uparrow}\{[0, 1]\setminus\{x\}: x\in [0, 1/2]\}$,
\end{center}
which has  no infimum in $\mathcal L$.

 Hence, by Lemma \ref{inf},  $\mathcal L$ fails to be a $\mathcal Q$-algebra. So $\mathcal Q$ cannot  give the free deflationary semilattice for $\mathbf{WF}$.
\end{exa}

\section{The commutativity between \texorpdfstring{$\mathcal H$}{H} and \texorpdfstring{$\mathcal Q$}{Q}}

\begin{defi}\label{phi}
 Let $L$ be a well-filtered dcpo. We define $\phi_{L}: \mathcal H (\mathcal Q(L))\rightarrow\mathcal Q(\mathcal H(L))$ and $\psi_{L}: \mathcal Q(\mathcal H(L))\rightarrow\mathcal H (\mathcal Q(L))$ as
 \begin{center}
 	 $\phi_{L}(\mathcal A) = \{A\in \mathcal H(L): \forall K\in \mathcal A, A\cap K\neq\emptyset\}$, for any $\mathcal A\in \mathcal H (\mathcal Q(L))$,
\\
	$\psi_{L}(\mathcal K) = \{K\in \mathcal Q(L): \forall A\in \mathcal K, K\cap A\neq\emptyset\}$, for any $\mathcal K\in \mathcal Q (\mathcal H(L))$.
\end{center}
\end{defi}
\noindent
Notably, these two definitions are similar in form to those given in \cite{brecht} and \cite{heck1991}, where the authors give the definitions and discuss them for all topological spaces and $\mathcal U_S$-admitting dcpos, respectively.    For convenience, we will omit the subscripts when the object acted on is clear.

\begin{lem}\label{coherent}
 Every complete lattice endowed with the Scott topology is coherent.
 \begin{proof}
 	\cite[Corollary 3.2]{Xi}  guarantees the well-filteredness of a complete lattice $L$. Since for any $x, y\in L$, ${\uparrow}x\cap {\uparrow}y = {\uparrow}(x\vee y)$,  by  \cite[Corollary 3.2]{note}, we have that $L$ is coherent.
 	\qedhere
 \end{proof}
\end{lem}

\begin{lem}
	The functions $\phi$ and $\psi$ defined in Definition \ref{phi} are well-defined.
	\begin{proof}
		To make it clear, let us define $\phi': \mathcal Q(L)\rightarrow \mathcal Q(\mathcal H(L))$ and $\psi': \mathcal H(L)\rightarrow \mathcal H(\mathcal Q(L))$ as follows:
		\begin{center}
		 $\phi'(K) = \{A\in \mathcal H(L):  A\cap K\neq\emptyset\}$, 	for any $K\in \mathcal Q(L)$,
			\\
		$\psi'(A) = \{K\in \mathcal Q(L):  K\cap A\neq\emptyset\}$, 	for any $A\in \mathcal H(L)$.
		\end{center}
Then $\phi(\mathcal A) = \bigcap\{\phi'(K): K\in \mathcal A\}$ and $\psi(\mathcal K) = \bigcap\{\psi'(A): A\in \mathcal K\}$.
One can easily check that $\phi'(K) = \mathcal Q\eta_L(K) = {\uparrow}\eta_{L} (K)$ and $\psi'(A) = \mathcal H\theta_L(A) = cl(\theta_{L}(A))$, which implies that $\phi'$ and $\psi'$ are well-defined. We know the intersection of arbitrary closed subsets is still closed. So $\psi$ is well-defined. Now we focus on $\phi$.  Lemma \ref{coherent} reveals that
 $\Sigma\mathcal H(L)$ is coherent. Combining with the well-filteredness of $\Sigma\mathcal H(L)$, one can obtain that each intersection of a collection of  compact saturated subsets in $\mathcal H(L)$ is still compact saturated.  Thus $\phi(\mathcal A)\in \mathcal Q(\mathcal H(L))$, i.e., $\phi$ is well-defined.
		\qedhere
	\end{proof}
\end{lem}
\begin{prop}
Let $L$ be a well-filtered dcpo. The function $\phi: \mathcal{HQ}(L)\rightarrow \mathcal{QH}(L)$ defined in Definition \ref{phi} preserves all sups, hence is Scott-continuous.
	\begin{proof}
Assume that there are $\mathcal A, \mathcal B\in \mathcal {HQ}(L)$ with $\mathcal A\leq \mathcal B$, i.e., $\mathcal A\subseteq \mathcal B$. By the definition of $\phi$, one can easily see that $\phi(\mathcal B)\subseteq \phi(\mathcal A)$, equivalently, $\phi(\mathcal A)\leq \phi(\mathcal B)$. So $\phi$ is order-preserving.
	Let $\{\mathcal A_i: i\in I\}$ be an arbitrary  subset of $\mathcal{HQ}(L)$. 	To complete this proof, we need to verify that $\phi(\bigvee_{i\in I}\mathcal A_i) = \bigvee_{i\in I} \phi(\mathcal A_i)$. Since $\phi$ is order-preserving, for each $i\in I$, $\phi(\mathcal A_i)\leq\phi(\bigvee_{i\in I} \mathcal A_i)$. Then  $\bigvee_{i\in I}\phi(\mathcal A_i)\leq\phi(\bigvee_{i\in I} \mathcal A_i)$ holds. Now we consider the reverse, more precisely, to prove  $  \bigcap_{i\in I} \phi(\mathcal A_i)\subseteq\phi(cl(\bigcup_{i\in I}\mathcal A_i))$.
		Take any closed subset $A\in  \bigcap_{i\in I} \phi(\mathcal A_i)$. Suppose that $A\notin \phi(cl(\bigcup_{i\in I}\mathcal A_i))$, i.e., there exists a compact saturated subset $K_0\in cl(\bigcup_{i\in I}\mathcal A_i)$ such that $A\cap K_0 = \emptyset$. Then $K_0\in \Box(L\setminus A)$, where $\Box(L\setminus A)$ is a  Scott open subset of  $\mathcal Q(L)$ by the well-filteredness of $L$. So $\Box(L\setminus A)\cap (\bigcup_{i\in I}\mathcal A_i)\neq\emptyset$, which implies that there is an $i_0\in I$ such that $\Box(L\setminus A)\cap \mathcal A_{i_0}\neq\emptyset$. That is to say, we can find a  $K_1\in \mathcal A_{i_0}$ contained in $L\setminus A$. This immediately indicates that $K_1\cap A = \emptyset$; thus $A\notin \phi(\mathcal A_{i_0})$, which contradicts $A\in \bigcap_{i\in I}\phi(\mathcal A_i)$. Therefore, $A\in \phi(cl(\bigcup_{i\in I}\mathcal A_i))$. From the arbitrariness of $A$, we conclude that  $  \bigcap_{i\in I} \phi(\mathcal A_i)\subseteq\phi(cl(\bigcup_{i\in I}\mathcal A_i))$, that is, $\phi(\bigvee_{i\in I} \mathcal A_i) \leq\bigvee_{i\in I}\phi(\mathcal A_i)$ holds.
		\qedhere
	\end{proof}
\end{prop}

\begin{defi}
A poset $P$ ($T_0$ space $X$)  is said to have \emph{property $(\mathrm{KC})$} if for every compact saturated subset $\mathcal K$ of $\Gamma(P)$  under the Scott topology ($\mathcal K$ of $\Gamma(X)$ with respect to the lower Vietoris topology), any $U\in \sigma(P)$ ($U\in \mathcal O(X)$)  that meets all members of $\mathcal K$ contains a compact saturated subset $K$ that still meets all members of $\mathcal K$.
\end{defi}

\noindent
Recall that a space $X$ is \emph{consonant} if and only if for any Scott open subset  $\mathcal U \subseteq \mathcal O(X)$ and any $U\in \mathcal U$, there exists a $K\in \mathcal Q(X)$ such that $U\in \Phi(K)\subseteq\mathcal U$, where $\Phi(K) = \{V\in \mathcal O(X): K\subseteq V\}$. The concept was proposed initially by Dolecki et al. in \cite{conso}  to answer some questions related to topologies on the hyperspace of closed subsets of a topological space. We say a poset is \emph{consonant} if it endowed with the Scott topology is consonant.

\begin{lem}\label{kc}
 If a poset $P$ is consonant, then it satisfies the property $(\mathrm{KC})$.
 \begin{proof}
Let $\mathcal K$ be a compact saturated subset of $(\Gamma(P), \sigma(\Gamma(P)))$ and $U$ a Scott open subset of $P$  satisfying $U\cap A\neq\emptyset$ for all $A\in \mathcal K$. Set $\mathcal U = \bigcap\{\vartriangle\!\! A: A\in \mathcal K\}$, where $\vartriangle\!\! A = \{V\in \sigma(P): V\cap A\neq\emptyset\}$. We claim that $\mathcal U\in \sigma(\sigma(P))$.  Obviously, it is an upper set. Assume $\mathcal V$ is a directed subset of $\sigma(P)$ and $\sup \mathcal V = \bigcup \mathcal V\in \mathcal U$. Then $\bigcup \mathcal V\in \vartriangle\!\! A$ for all $A\in \mathcal K$, that is, $\bigcup \mathcal V\cap A\neq\emptyset$ for all $A\in \mathcal K$. It means that for each $A\in \mathcal K$, there exists a $V_{A}\in \mathcal V$ such that $V_{A}\cap A\neq\emptyset$, i.e., $A\in \lozenge V_{A}$, and thus $\mathcal K\subseteq \bigcup\{\lozenge V: V\in \mathcal V\}$.  By the compactness of $\mathcal K$ and the directedness of $\{\lozenge V: V\in \mathcal V\}$, there must exist a $V_{0}\in \mathcal V$ such that $\mathcal K\subseteq \lozenge V_{0}$. It follows that $V_{0}\in \vartriangle\!\! A$ for each $A\in \mathcal K$, that is to say, $V_{0}\in \mathcal U$. Thus $\mathcal U$ is Scott open. Since $P$ is consonant and $U\in\mathcal U$, there exists a $K\in \mathcal Q(P)$ satisfying $U\in \Phi(K)\subseteq \mathcal U$. We will complete the proof if we show that $K\cap A\neq\emptyset$ for all $A\in \mathcal K$. Suppose $K\cap A_{0} = \emptyset$ for some $A_{0}\in \mathcal K$. Then $K\subseteq P\setminus A_{0}$, where $ P\setminus A_{0}\in \sigma(P)$. This implies   $P\setminus A_{0}\in \Phi(K)\subseteq \mathcal U$. So $P\setminus A_{0}\in \;\vartriangle\!\! A_{0}$,  which contradicts the definition of $\vartriangle\!\! A_{0}$.
 	\qedhere
 \end{proof}
\end{lem}

\noindent
Note that the property $(\mathrm{KC})$ for $T_0$ spaces  is exactly the condition
	``$\mathcal K \in\Box\diamondsuit U$ implies $\tau (\mathcal K) \in \diamondsuit \Box U$'', i.e., $\Box\diamondsuit U \subseteq  \tau^{-1}(\diamondsuit\Box U)$
	in the notation of  \cite{brecht}, where the $\tau$ in \cite{brecht} is the same as the $\psi$ defined in
	Definition \ref{phi}. Since Lemma 6.5 in \cite{brecht} has shown that $\tau^{-1}(\diamondsuit\Box U) \subseteq   \Box\diamondsuit U$, the
	following lemma  is equivalent to the fact that (1) $\Leftrightarrow$ (3) in
	\cite[Theorem 6.13]{brecht}. 
\begin{lem}
	A $T_{0}$ space $X$ is consonant if and only if it has the property $(\mathrm{KC})$.

\end{lem}

\noindent
Now see  \cite[Theorem 6.10]{brecht}, it states that for a space $X$, $\sigma(\mathcal O(X))$ and $\mathcal Q(\Gamma(X))$ are isomorphic lattices. Recently, Chen et al. in \cite{chenyu} proved that if a poset $P$ endowed with the Scott topology is core-compact, then $\sigma(\Gamma(P)) = \nu(\Gamma(P))$. Thus for a poset $P$, $\sigma(\sigma(P))$ would be isomorphic to $\mathcal Q(\Gamma(P))$ as long as $\Sigma P$ is core-compact. Combining all of these results and the above lemma, we obtain the following result.

\begin{cor}
	Let $P$ be a poset. If $(P, \sigma(P))$ is core-compact and has the property $(\mathrm{KC})$, then $P$ is consonant.
\end{cor}

\begin{thm}\label{iso}
	Let $L$ be  a well-filtered dcpo. Then $\mathcal H(\mathcal Q(L))$ is isomorphic to $\mathcal Q(\mathcal H(L))$ under the maps $\phi$ and $\psi$ if and only if the following conditions are satisfied.
	\begin{enumerate}[(i)]
		\item $\sigma(Q(L))\subseteq \upsilon(Q(L));$ and
		\item $L$ has the property $(\mathrm{KC})$.
	\end{enumerate}
\begin{proof}
	\leavevmode
	$(\Rightarrow)$: 
	(i) Let $\mathcal U$ be a Scott open subset of $\mathcal Q(L)$ and $K$ an element in $\mathcal U$. Then $K\notin \mathcal Q(L)\setminus \mathcal U$, where $\mathcal Q(L)\setminus \mathcal U$ is a  closed subset of $\Sigma\mathcal Q(L)$. Since $\psi$ is surjective, there exists a $\mathcal K\in \mathcal Q(\mathcal H(L))$ such that $\mathcal Q(L)\setminus \mathcal U = \psi(\mathcal K)$; consequently, $K\notin \psi(\mathcal K)$. From the definition of $\psi$, we have $K\cap A = \emptyset$ for some $A\in \mathcal K$. Set $V = L \setminus A$. Then $V\in \sigma(L)$ and   $K\subseteq V$, i.e., $K\in \Box V$. We now show  $\Box V\subseteq\mathcal U$. Take any $M\in \Box V$. Then $M\cap A = \emptyset$, which implies $M\notin \psi(\mathcal K)  = \mathcal Q(L)\setminus \mathcal U$;  hence, $M\in \mathcal U$. Therefore $K\in \Box V\subseteq\mathcal U$, so $\mathcal U\in \upsilon(Q(L))$. By the arbitrariness of $\mathcal U$, we conclude that $\sigma(Q(L))\subseteq \upsilon(Q(L))$.

 (ii) Let $\mathcal K$ be a Scott compact saturated subset of $\Gamma(L)$ and $U\in \sigma(L)$ with $U\cap A\neq\emptyset$ for all $A\in \mathcal K$. Since $U \cap (L \setminus U) = \emptyset$, we have $L \setminus U \notin\mathcal K$.
	By the surjectivity of $\phi$, there exists an $\mathcal A\in \mathcal H(\mathcal Q(L))$ such that $\mathcal K = \phi(\mathcal A)$. Thus $L\setminus U\notin \phi (\mathcal A)$. This means $(L\setminus U)\cap K_{0} = \emptyset$, i.e., $K_{0}\subseteq U$ for some $K_{0}\in \mathcal A$. Now suppose that $K_{0}\cap A_{1} = \emptyset$ for some $A_{1}\in \mathcal K$. Then $A_{1}\notin\phi (\mathcal A) = \mathcal K$, which is a contradiction  obviously. Therefore, $K_{0}$ intersects every element in $\mathcal K$, and $L$ has the property $(\mathrm{KC})$ as desired.
	
	$(\Leftarrow)$: Clearly, $\phi$ and $\psi$ are order-preserving. 	Since order-preserving isomorphisms between dcpos are Scott-continuous, the proof will be complete
once we show that $\psi\circ\phi$ and  $\phi\circ\psi$ are identity maps, more precisely, the following two equations hold for arbitrary $\mathcal A\in \mathcal H(\mathcal Q(L))$ and $\mathcal K\in \mathcal Q(\mathcal H(L))$:
	\begin{center}
		$\psi\circ\phi(\mathcal A) = \{K\in \mathcal Q(L): \forall A\in \phi(\mathcal A), K\cap A\neq\emptyset\} = \mathcal A$,
		\\
	$\phi\circ\psi(\mathcal K) = \{A\in \mathcal H(L): \forall K\in \psi(\mathcal K), A\cap K\neq\emptyset\} = \mathcal K$.	
	\end{center}
It follows directly from the definitions of $\phi$ and $\psi$  that $\mathcal A\subseteq\psi\circ\phi(\mathcal A)$ and $\mathcal K\subseteq\phi\circ\psi(\mathcal K)$. Thus we just need to check the reverse inclusions.
	
	For the first inclusion, suppose that there is a compact saturated subset $K$  belonging to $\psi\circ\phi(\mathcal A)$ but not to $\mathcal A$. Then $K\in \mathcal Q(L)\setminus \mathcal A$, which is  Scott open  in $\mathcal Q(L)$ and by $(\mathrm i)$, is also  open in the upper Vietoris topology. Hence, there exists a $U\in \sigma(L)$ such that $K\in \Box U\subseteq\mathcal Q(L)\setminus \mathcal A$, which implies $K\cap (L \setminus U) = \emptyset$. This would contradict the fact that $K\in \psi\circ\phi(\mathcal A)$ once we  show $L\setminus U\in \phi(\mathcal A)$. Assume $(L\setminus U)\cap K_{0} = \emptyset$ for some $K_{0}\in \mathcal A$. Then $K_{0}\in \Box U$, from which we have $K_{0}\notin \mathcal A$, a contradiction. Thus $L\setminus U\in \phi(\mathcal A)$;  consequently, $\psi\circ\phi(\mathcal A)\subseteq \mathcal A$ holds.
	
To prove the second inclusion, for the sake of a contradiction, we  assume  that there is an $A\in \phi\circ\psi(\mathcal K)$ but $A\notin \mathcal K$. This indicates that there is no $C\in \mathcal K$  contained in $A$, that is, $(L\setminus A)\cap C\neq\emptyset$ for all $C\in \mathcal K$. Since $L$ has the property $(\mathrm{KC})$, we can find a $K_{0}\subseteq L\setminus A$ that intersects each $C$ in $\mathcal K$. Thus $K_{0}\in \psi(\mathcal K)$ and $K_{0}\cap A = \emptyset$, which contradicts that $A\in \phi\circ\psi(\mathcal K)$. So each element in $\phi\circ\psi(\mathcal K)$ belongs to $\mathcal K$, that is, $\phi\circ\psi(\mathcal K)\subseteq\mathcal K$.
	\qedhere
\end{proof}

\end{thm}
\noindent
It was A. Schalk who proved that for a locally compact sober space $X$, the upper Vietoris topology and the Scott topology  on $\mathcal Q(X)$ coincide (see \cite[Lemma 7.26]{as}), meanwhile, \cite[Proposition 5.4]{chenyu} reveals that each locally compact space is consonant. Therefore, combining Lemma \ref{kc} and Theorem \ref{iso},  we further  draw the conclusion in the  following.
\begin{cor}\label{lsi}
	If $L$ is a locally compact sober dcpo, then $\mathcal {QH}(L)$ is isomorphic to $\mathcal {HQ}(L)$ under the function $\phi$ (or $\psi$).
\end{cor}
\section{The double power well-filtered dcpo}
In the study of monads, a question always worth considering is the following: given two monads with underlying functors $\mathcal  S$ and $\mathcal T$ on a category $\mathbf C$, does $ \mathcal T\circ \mathcal S$ (or $\mathcal S\circ \mathcal T$) also carry  a monad structure? In 1969, J. Beck introduced the concept of \emph{distributive laws} and showed that their existence  provides a sufficient condition for such a composition  to form a monad \cite{beck}.

\begin{defiC}[\cite{beck}]\label{dl}
	Let $\mathcal S = (\mathcal S, \eta^{\mathcal S}, \mu^{\mathcal S})$ and $\mathcal T = (\mathcal T, \eta^{\mathcal T}, \mu^{\mathcal T})$ be two monads on a category $\mathbf C$. A \emph{distributive law} of $\mathcal S$ over $\mathcal T$ is a natural transformation $ l: \mathcal{TS}\rightarrow \mathcal{ST} $ such that the diagrams in Figure 1 commute.	
\begin{figure}[ht]
	\centering
	\begin{tikzpicture}[line width=0.6pt,scale=0.7]
	\node[right]  at(-5,0) (a) { $\mathcal T$};
    \node[right]  at(-7,-2) (b) { $\mathcal {TS}$};
    \draw[->](a)--(b);
    \node[right]  at(-6.4,-0.75)  { $\mathcal T\eta^{\mathcal S}$};
    \node[right]  at(-3,-2) (c) { $\mathcal{ST}$};
    \draw[->](a)--(c);
    \draw[->](b)--(c);
    \node[right]  at(-3.7,-0.75)  { $\eta^{\mathcal S}_{\mathcal T}$};
    \node[right]  at(-4.8,-1.7)  { $l$};

    \node[right]  at(5,0)  (d) { $\mathcal S$};
    \node[right] at(3,-2) (e) { $\mathcal{TS}$};
    \node[right] at(7,-2)  (f){$\mathcal{ST}$};
    \draw[->](d)--(e);
    \draw[->](d)--(f);
    \draw[->](e)--(f);
     \node[right]  at(6.3,-0.75)  {$\mathcal S\eta^{\mathcal T}$};
      \node[right] at(3.7,-0.75)  {$\eta^{\mathcal T}_{\mathcal S}$};
       \node[right]  at(5.2,-1.7)  { $l$};
  \node[right]  at(-5,-4) (h) { $\mathcal{TSS}$};
  \node[right]  at(0,-4) (i) { $\mathcal {STS}$};
   \node[right] at(5,-4) (j) { $\mathcal {SST}$};
     \draw[->](h)--(i);
     \draw[->](i)--(j);
    \node[right] at(-2.4,-3.7)  {$l_{\mathcal S}$};
    \node[right]  at(2.6,-3.7)  { $\mathcal Sl$};
 \node[right]  at(-4.8,-6.3) (k) { $\mathcal{TS}$};
  \node[right]  at(5.2,-6.3) (l) { $\mathcal{ST}$};
 \draw[->](h)--(k);
 \draw[->](j)--(l);
 \draw[->](k)--(l);
 \node[right]  at(-5.5,-5.2)  { $T\mu^{\mathcal S}$};
  \node[right]  at(5.8,-5.2)  { $\mu^{\mathcal S}_{\mathcal T}$};
  \node[right]  at(0.4,-6)  { $l$};

  \node[right]  at(-5,-8) (o) {$\mathcal{TTS}$};
  \node[right] at(0,-8) (p) { $\mathcal{TST}$};
  \node[right] at(5,-8) (q) { $\mathcal{STT}$};
  \draw[->](o)--(p);
  \draw[->](p)--(q);
  \node[right]  at(-2.4,-7.7)  { $\mathcal Tl$};
  \node[right] at(2.6,-7.7)  {$l_{\mathcal T}$};
  \node[right] at(-4.8,-10.3) (r) { $\mathcal{TS}$};
  \node[right]  at(5.2,-10.3) (s) {$\mathcal{ST}$};
  \draw[->](o)--(r);
  \draw[->](q)--(s);
  \draw[->](r)--(s);
  \node[right]  at(-5.3,-9.2)  { $\mu^{\mathcal T}_{\mathcal S}$};
  \node[right]  at(5.8,-9.2)  { $\mathcal S\mu^{\mathcal T}$};
  \node[right] at(0.4,-10)  {$l$};
	\end{tikzpicture}	
		\scriptsize    \\ \small{Figure 1}
\end{figure}
\end{defiC}

\begin{defiC}[\cite{ct}]
	Let $\mathcal T = (\mathcal T, \eta, \mu)$ and $\mathcal T' = (\mathcal T', \eta', \mu')$ be two monads on a category $\mathbf C$. A \emph{morphism of monads} is a natural transformation $\lambda: \mathcal T\rightarrow \mathcal T'$ such that the   diagrams in Figure 2  commute.
	\begin{figure}[ht]
		\centering
		\begin{tikzpicture}[line width=0.6pt, scale=0.7]
			\node[right]  at(-2.5,0)(a)  { $Id$};
			\node[right]  at(1,1.1) (b)  { $\mathcal T$};
			\node[right]  at(1,-1.1) (c){ $\mathcal T'$};
			\draw[->](a)--(b);
			\draw[->](a)--(c);
			\draw[->](b)--(c);
			\node[right]  at(-0.7,0.8)  { $\eta$};
			\node[right]  at(-0.7,-0.9)  { $\eta'$};
			\node[right] at(1.3,0)  {$\lambda$};
			
			\node[right]  at(4,1.1)(d)  {$\mathcal {TT}$};
			\node[right]  at (7,1.1) (e){ $\mathcal{T'T}$};	
			\node[right]  at(10,1.1) (f)  {$\mathcal {T'T'}$};
			\node[right]  at(4.2,-1.1) (g)  { $\mathcal T$};
			\node[right]  at(10.3,-1.1) (h)  { $\mathcal T'$};
			\draw[->](d)--(e);
			\draw[->](e)--(f);
			\draw[->](d)--(g);
			\draw[->](f)--(h);
			\draw[->](g)--(h);
			\node[right] at(5.5,1.4)  { $\lambda_{\mathcal T}$};
			\node[right]  at(8.5,1.4)  { $\mathcal T'\lambda$};
			\node[right] at(4,0)  { $\mu$};
			\node[right]  at(10.7,0)  { $\mu'$};
			\node[right] at(7.2,-0.8)  { $\lambda$};
			
		\end{tikzpicture}
		\vspace{0.45em}
				\scriptsize    \\ \small{Figure 2}
	\end{figure}
\end{defiC}

\noindent
Particularly, the morphisms of monads are called \emph{triple maps} by J. Beck in \cite{beck}. Meanwhile, J. Beck also proved the following theorem.
\begin{thmC}[\cite{beck}]\label{beck1}
	Let $\mathcal S = (\mathcal S, \eta^{\mathcal S}, \mu^{\mathcal S})$, $\mathcal T = (\mathcal T, \eta^{\mathcal T}, \mu^{\mathcal T})$ be two monads on a category $\mathbf C$. Then the following statements are equivalent:
	\begin{enumerate}[(1)]
	\item There are distributive laws $d: \mathcal{ST}\rightarrow \mathcal {TS};$
	
	\item	There are multiplications $\mu: \mathcal{TSTS}\rightarrow \mathcal{TS}$ such that
	\begin{itemize}
		\item $(\mathcal{TS}, \eta_{\mathcal S}^{\mathcal T}\eta^{\mathcal S}, \mu)$ is a monad;
		\item The natural transformations $ \eta^{\mathcal T}_{\mathcal S}: \mathcal S\rightarrow \mathcal {TS}$ and $\mathcal T\eta^{\mathcal S}: \mathcal T\rightarrow \mathcal{TS}$ are morphisms of monads.
	\item The middle unitary law 
	\begin{figure}[ht]
		\centering
		\begin{tikzpicture}[line width=0.6pt,scale=0.95]
			\node[right]  at(-5,0) (a) { $\mathcal {TS}$};
			\node[right]  at(-7,-2) (b) { $\mathcal {TSTS}$};
			\draw[->](a)--(b);
			\node[right]  at(-7.7,-0.8)  { $\mathcal T\eta^{\mathcal S}_{\mathcal {TS}}\circ\mathcal T\eta^{\mathcal T}_{\mathcal S}$};
			\node[right]  at(-3,-2) (c) { $\mathcal{TS}$};
			\draw[->](a)--(c);
			\draw[->](b)--(c);
			\node[right]  at(-3.7,-0.8)  { = };
			\node[right]  at(-4.6,-1.8)  { $\mu$};
		\end{tikzpicture}
	\end{figure}
	\\
i.e., $\mu\circ\mathcal T\eta^{\mathcal S}_{\mathcal {TS}}\circ\mathcal T\eta^{\mathcal T}_{\mathcal S} = id_{\mathcal{TS}}$ holds.
	\end{itemize}

		\end{enumerate}
\end{thmC}
\noindent
From the proposition presented in Section 2 of J. Beck \cite{beck}, a more general result concerning the morphisms of monads  was  established by T. Fritz et al.

	\begin{lemC}[\cite{tps}]\label{algebra}
		Let $\mathcal S = (\mathcal S, \eta_{\mathcal S}, \mu_{\mathcal S})$ and $\mathcal T = (\mathcal T, \eta_{\mathcal T}, \mu_{\mathcal T})$ be two monads on a category $\mathbf C$ and $\lambda: \mathcal S\rightarrow \mathcal T$ a morphism of monads. Then every $\mathcal T$-algebra $(A, \alpha)$ can be equipped with an $\mathcal S$-algebra structure via $(A, \alpha)\mapsto (A, \alpha\circ\lambda)$. Moreover, a $\mathcal T$-algebra homomorphism $f: (A, \alpha)\rightarrow (B, \beta)$ induces an $\mathcal S$-algebra homomorphism  $f: (A, \alpha\circ\lambda)\rightarrow (B, \beta\circ\lambda)$ in a functorial way.
	\end{lemC}
\begin{prop}\label{dis}
$\phi$ : $\mathcal H\mathcal Q \rightarrow \mathcal Q\mathcal H$  defined in Definition \ref{phi} is a distributive law of $\mathcal Q$ over $\mathcal H$.
\begin{proof}
$\bf {Claim ~1:}$ $\phi$ is a natural transformation.
	
		Let $L$ and $M$ be well-filtered dcpos and  $f: L\rightarrow M$  a Scott-continuous map. We need to show that the following diagram commutes:
		\vspace{4em}
		\begin{figure}[ht]
			\centering
			\begin{tikzpicture}[line width=0.6pt,scale=0.7]

				\node[right]  at(-5,0) (o) {$\mathcal{HQ}(L)$};
				\node[right] at(0,0) (p) { $\mathcal{QH}(L)$};
			
				\draw[->](o)--(p);
				\node[right]  at(-2,0.3)  { $\phi_L$};
				
				\node[right] at(-5,-3.2) (r) { $\mathcal{HQ}(M)$};
				\node[right]  at(0,-3.2) (s) {$\mathcal{QH}(M)$};
				\draw[->](o)--(r);
				
				\draw[->](r)--(s);
				\draw[->](p)--(s);
				\node[right]  at(-5.8,-1.6)  { $\mathcal {HQ}(f)$};
				\node[right]  at(0.9,-1.6)  { $\mathcal{QH}(f)$};
				\node[right] at(-2,-2.9)  {$\phi_M$};
			\end{tikzpicture}
		\end{figure}
		\\
	i.e., for any $\mathcal A\in \mathcal {HQ}(L)$, $\mathcal {QH}(f)\circ\phi_L(\mathcal A) = 	\phi_M\circ\mathcal {HQ}(f)(\mathcal A) $, where 
	\begin{center}
	$\mathcal {QH}(f)\circ\phi_L(\mathcal A) = {\uparrow }\{cl (f(A)): \forall K\in \mathcal A, A\cap K\neq\emptyset\}$,
	
	$\phi_M\circ\mathcal {HQ}(f)(\mathcal A) = \{B\in \mathcal H(M): \forall N\in cl\{{\uparrow}f(K): K\in \mathcal A\}, B\cap N\neq\emptyset\}$.
	\end{center}
First, suppose that there exists a $cl(f(A))
\in \mathcal {QH}(f)\circ\phi_L(\mathcal A)$ but $cl(f(A))
\notin \phi_M\circ\mathcal {HQ}(f)(\mathcal A)$. Then we can find an $N_0\in cl\{{\uparrow}f(K): K\in \mathcal A\}$ satisfying  $cl(f(A))\cap N_0 = \emptyset$. It follows that $N_0\in \Box(M\setminus cl(f(A))) $. Since the well-filteredness of $M$ implies that  $\Box(M\setminus cl(f(A))) \in \sigma(\mathcal Q(M))$,  $\Box(M\setminus cl(f(A)))\cap \{{\uparrow}f(K): K\in \mathcal A\}\neq\emptyset$. So there is  a $K_0\in \mathcal A$ such that ${\uparrow }f(K_0)\subseteq M\setminus cl(f(A))$, that is to say, ${\uparrow }f(K_0)\cap cl(f(A)) =  \emptyset$. Then $K_0\cap A = \emptyset$, which obviously contradicts  the fact that $A$ meets every element of $\mathcal A$. Thus, $\{cl (f(A)): \forall K\in \mathcal A, A\cap K\neq\emptyset\}\subseteq \phi_M\circ\mathcal {HQ}(f)(\mathcal A)$;  therefore, $\mathcal {QH}(f)\circ\phi_L(\mathcal A) \subseteq 	\phi_M\circ\mathcal {HQ}(f)(\mathcal A) $.

Now we prove the reverse inclusion.  Let $B\in \phi_M\circ\mathcal {HQ}(f)(\mathcal A)$. Then for each $K\in\mathcal A$, $B\cap {\uparrow}f(K)\neq\emptyset$. Because $B$ is a down set, this implies $B\cap f(K)\neq\emptyset$ for all $K\in\mathcal A$. We pick $x_K\in K$ such that $f(x_K)\in B\cap f(K)$ and set $A_0 = cl(\{x_K: K\in \mathcal A\})$. Then $cl(f(A_0)) =  cl(f(\{x_K: K\in \mathcal A\})) = cl(\{f(x_K): K\in \mathcal A\})\subseteq B$ using that $B$ is closed. Meanwhile, $A_0\cap K\neq\emptyset$ for each $K\in \mathcal A$. So $B\in 	\mathcal {QH}(f)\circ\phi_L(\mathcal A)$; thus, $\phi_M\circ\mathcal {HQ}(f)(\mathcal A)\subseteq\mathcal {QH}(f)\circ\phi_L(\mathcal A)$ holds.

$\bf {Claim ~2:}$ The first diagram in Definition \ref{dl} commutes.

Concretely, given a well-filtered dcpo $L$, we need to demonstrate that the following diagram commutes.
\begin{figure}[ht]
	\centering
	\begin{tikzpicture}[line width=0.6pt,scale=0.7]
		\node[right]  at(-5,0) (a) { $\mathcal H(L)$};
		\node[right]  at(-7,-2) (b) { $\mathcal {HQ}(L)$};
		\draw[->](a)--(b);
		\node[right]  at(-6.2,-0.75)  { $\mathcal H\theta_L$};
		\node[right]  at(-3,-2) (c) { $\mathcal{QH}(L)$};
		\draw[->](a)--(c);
		\draw[->](b)--(c);
		\node[right]  at(-3.3,-0.75)  { $\theta_{\mathcal H(L)}$};
		\node[right]  at(-4.4,-1.7)  { $\phi_L$};
	\end{tikzpicture}
\end{figure}
\\
Let $A\in \mathcal H(L)$. Then 
$\theta_{\mathcal H(L)}(A) = {\uparrow}\{A\}$ and
\begin{center}
$\phi_L\circ\mathcal H\theta_L(A) = \{B\in \mathcal H(L): \forall K\in cl(\{{\uparrow}x: x\in A\}): B\cap K\neq\emptyset\}$.
\end{center}
First, suppose for contradiction that there exists a $K\in cl(\{{\uparrow}x: x\in A\})$ such that $A\cap K = \emptyset$. Then $K\subseteq L\setminus A$, i.e., $K\in \Box (L\setminus A)$. By the well-filteredness of $L$, we know $\Box (L\setminus A)\in \sigma(\mathcal Q(L))$. So $\Box (L\setminus A)\cap \{{\uparrow}x: x\in A\}\neq\emptyset$, that is to say, we can find an $x\in A$ such that ${\uparrow}x\subseteq L\setminus A$, which implies $x\notin A$,  a contradiction. So for every $K\in cl(\{{\uparrow}x: x\in A\}), A\cap K\neq\emptyset$. From the definition of $\phi_L$, we have  $A\in \phi_L\circ\mathcal H\theta_L(A)$;  hence, $\theta_{\mathcal H(L)}(A) = {\uparrow}\{A\}\subseteq \phi_L\circ\mathcal H\theta_L(A)$.

Now we consider the reverse inclusion. Take any $B\in \phi_L\circ\mathcal H\theta_L(A)$. Then for every $x\in A$, $B\cap {\uparrow}x\neq\emptyset$. It means that each $x\in A$  belongs to $B$. Thus $A\subseteq B$,  equivalently, $B\in {\uparrow}\{A\}$. By the arbitrariness of $B$, we conclude that $\phi_L\circ\mathcal H\theta_L(A)\subseteq \theta_{\mathcal H(L)}(A)$.

$\bf {Claim ~3:}$ The second diagram in Definition \ref{dl} commutes.

Specifically, given a well-filtered dcpo $L$, we need to prove that the following diagram commutes, 
\begin{figure}[ht]
	\centering
	\begin{tikzpicture}[line width=0.6pt,scale=0.7]
		\node[right]  at(-5,0) (a) { $\mathcal Q(L)$};
		\node[right]  at(-7,-2) (b) { $\mathcal {HQ}(L)$};
		\draw[->](a)--(b);
		\node[right]  at(-6.3,-0.75)  { $\eta_{\mathcal Q(L)}$};
		\node[right]  at(-3,-2) (c) { $\mathcal{QH}(L)$};
		\draw[->](a)--(c);
		\draw[->](b)--(c);
		\node[right]  at(-3.3,-0.75)  { $\mathcal Q\eta_L$};
		\node[right]  at(-4.4,-1.7)  { $\phi_L$};
	\end{tikzpicture}
\end{figure}
\\
i.e., to verify that $\phi_L\circ\eta_{\mathcal Q(L)}(K) = \mathcal Q\eta_L(K)$ for any $K\in \mathcal Q(L)$, where $\mathcal Q\eta_L(K) = {\uparrow}\{{\downarrow}x: x\in K\}$ and 
\begin{center}
$\phi_L\circ\eta_{\mathcal Q(L)}(K) = \{A\in \mathcal H(L): \forall M\in {\downarrow}\{K\}: A\cap M\neq\emptyset\}$.
\end{center}
Let  $A$ be an arbitrary element in $\phi_L\circ\eta_{\mathcal Q(L)}(K)$. Then $A\cap K\neq\emptyset$;  thus, there is an $x\in A\cap K$. By the closedness of $A$, we have ${\downarrow}x\subseteq A$, i.e., $A\in {\uparrow}\{{\downarrow}x\}$. So $A\in \mathcal Q\eta_L(K)$. The arbitrariness of $A$ guarantees that  $\phi_L\circ\eta_{\mathcal Q(L)}(K)\subseteq \mathcal Q\eta_L(K)$.

Now we focus on the reverse inclusion. For each $x\in K$, we have ${\downarrow}x\cap K\neq\emptyset$. So for any $M\in  {\downarrow}\{K\}$, i.e., $K\subseteq M$, ${\downarrow}x\cap M\neq \emptyset$ always holds. This means $\{{\downarrow}x: x\in K\}\subseteq \phi_L\circ\eta_{\mathcal Q(L)}(K) $;   thus, $\mathcal Q\eta_L(K)\subseteq \phi_L\circ\eta_{\mathcal Q(L)}(K)$.

$\bf {Claim ~4:}$ The third diagram in Definition \ref{dl} commutes, i.e., the following one commutes for a well-filtered dcpo $L$:
	\begin{figure}[ht]
		\centering
		\begin{tikzpicture}[line width=0.6pt,scale=0.7]
		
			\node[right]  at(-5,-4) (h) { $\mathcal{HQQ}(L)$};
			\node[right]  at(0,-4) (i) { $\mathcal {QHQ}(L)$};
			\node[right] at(5,-4) (j) { $\mathcal {QQH}(L)$};
			\draw[->](h)--(i);
			\draw[->](i)--(j);
			\node[right] at(-2.4,-3.7)  {$\phi_{\mathcal Q(L)}$};
			\node[right]  at(2.6,-3.7)  { $\mathcal Q\phi_L$};
			\node[right]  at(-4.8,-6.3) (k) { $\mathcal{HQ}(L)$};
			\node[right]  at(5.2,-6.3) (l) { $\mathcal{QH}(L)$};
			\draw[->](h)--(k);
			\draw[->](j)--(l);
			\draw[->](k)--(l);
			\node[right]  at(-5,-5.2)  { $\mathcal H\iota_L$};
			\node[right]  at(6.1,-5.2)  { $\iota_{\mathcal H(L)}$};
			\node[right]  at(0.4,-6)  { $\phi_L$};
		\end{tikzpicture}
	\end{figure}
	\\
For any $\mathbf A\in \mathcal{HQQ}(L)$, by computation, we have
\begin{center}
	
	$\iota_{\mathcal H(L)}\circ\mathcal Q\phi_L\circ\phi_{\mathcal Q(L)}(\mathbf A) = \bigcup{\uparrow}\{\phi_L(\mathcal A): \forall \mathcal K\in\mathbf A, \mathcal A\cap \mathcal K\neq\emptyset\}$,
	
$\phi_L\circ \mathcal H\iota_L(\mathbf A) = \{A\in \mathcal H(L): \forall K\in cl(\{\bigcup\mathcal K: \mathcal K\in \mathbf A\}), A\cap K \neq\emptyset\}$.
\end{center}
Consider each $\phi_L(\mathcal A) = \{B\in \mathcal H(L): \forall M\in \mathcal A, B\cap M\neq\emptyset\}$, in which $\mathcal A\in \mathcal{HQ}(L)$ and meets every element of $\mathbf A$. Suppose that there exists a $B\in\mathcal H(L)$  belonging to $\phi_L(\mathcal A)$ but not  to $\phi_L\circ \mathcal H\iota_L(\mathbf A)$, i.e., there is a $K\in cl(\{\bigcup\mathcal K: \mathcal K\in \mathbf A\})$ satisfying  $B\cap K = \emptyset$. Then $K\in \Box(L\setminus B)$. Applying the well-filteredness of $L$, we obtain that $\Box(L\setminus B)$ is a Scott open subset of $\mathcal Q(L)$;  further, we have $\Box(L\setminus B)\cap \{\bigcup\mathcal K: \mathcal K\in \mathbf A\}\neq\emptyset$. So we can find a $\mathcal K\in\mathbf A$ satisfying  $\bigcup\mathcal K\subseteq L\setminus B$,  which indicates that for any $C\in \mathcal K$, $C\cap B = \emptyset$. Combining with the fact that $B\in \phi_L(\mathcal A)$, one can conclude that there is no element of $\mathcal K$ belonging to $\mathcal A$, that is to say, $\mathcal K\cap\mathcal A = \emptyset$. This contradicts the assumption that $\mathcal A$ meets every element of $\mathbf A$. Hence, each $\phi_L(\mathcal A)$ is contained in 	$\phi_L\circ \mathcal H\iota_L(\mathbf A)$; accordingly, $\iota_{\mathcal H(L)}\circ\mathcal Q\phi_L\circ\phi_{\mathcal Q(L)}(\mathbf A)\subseteq \phi_L\circ \mathcal H\iota_L(\mathbf A)$.

Take an arbitrary element $A\in \phi_L\circ \mathcal H\iota_L(\mathbf A)$. Then for each $\mathcal K\in \mathbf A, A\cap \bigcup\mathcal K\neq\emptyset$; hence, we can choose an element $a_{\mathcal K}\in A\cap \bigcup\mathcal K$. This implies  that there exists some $K\in \mathcal K$ such that ${\uparrow}a_{\mathcal K}\subseteq K$, i.e., $K\leq{\uparrow}a_{\mathcal K}$. Since $\mathcal K$ is an upper set, ${\uparrow}a_{\mathcal K}\in \mathcal K, \; \forall \mathcal K\in \mathbf A$. Set $C = cl\{a_{\mathcal K}: \mathcal K\in \mathbf A\}$. Then  the mapping $\mathcal H\theta_L$ sends the set $C$ to  a Scott closed subset $cl(\{{\uparrow}c: c\in C\})$ of $\mathcal Q(L)$.  We now verify that $A\in \phi_L(cl(\{{\uparrow}c: c\in C\}))$ and $\phi_L(cl(\{{\uparrow}c: c\in C\}))\in \mathcal Q\phi_L\circ\phi_{\mathcal Q(L)}(\mathbf A)$. First, suppose that there exists an $M\in cl(\{{\uparrow}c: c\in C\})$ satisfying $A\cap M = \emptyset$. Then $M\in \Box(L\setminus A)$. Since $L$ is well-filtered, $\Box(L\setminus A)$ is a Scott open subset of $\mathcal Q(L)$, which implies $\Box(L\setminus A)\cap \{{\uparrow}c: c\in C\}\neq\emptyset$. In other words, there exists a $c\in C$ such that ${\uparrow}c\subseteq L\setminus A$;  hence,  $c\in L\setminus A$. Note that $L\setminus A\in \sigma(L)$ and $c\in cl\{a_{\mathcal K}: \mathcal K\in \mathbf A\}$, so we have $(L\setminus A)\cap \{a_{\mathcal K}: \mathcal K\in \mathbf A\}\neq\emptyset$.
Thus we can find an $a_{\mathcal K}$  that belongs to $ L\setminus A$, i.e., $a_{\mathcal K}\notin A$, which obviously leads to a contradiction. Therefore, $A\in \phi_L(cl(\{{\uparrow}c: c\in C\}))$.  Next,  for each $a_{\mathcal K}\in C$, we have ${\uparrow}a_{\mathcal K}\in cl(\{{\uparrow}c: c\in C\})$ and ${\uparrow}a_{\mathcal K}\in\mathcal K$ for every $\mathcal K\in \mathbf A$. It follows that $cl(\{{\uparrow}c: c\in C\})\cap\mathcal K\neq\emptyset, \;\forall \mathcal K\in \mathbf A$. This immediately yields $\phi_L(cl(\{{\uparrow}c: c\in C\}))\in \mathcal Q\phi_L\circ\phi_{\mathcal Q(L)}(\mathbf A)$. In summary, we find an element $\phi_L(cl(\{{\uparrow}c: c\in C\}))$ in $\mathcal Q\phi_L\circ\phi_{\mathcal Q(L)}(\mathbf A)$  such that   $A$  belongs to it. By the arbitrariness of $A$,  we conclude that $\phi_L\circ \mathcal H\iota_L(\mathbf A)\subseteq \iota_{\mathcal H(L)}\circ\mathcal Q\phi_L\circ\phi_{\mathcal Q(L)}(\mathbf A)$.

$\bf {Claim ~5:}$ The fourth diagram in Definition \ref{dl} commutes, i.e., the following one commutes:
\begin{figure}[ht]
	\centering
	\begin{tikzpicture}[line width=0.6pt,scale=0.7]
		
		\node[right]  at(-5,-4) (h) { $\mathcal{HHQ}(L)$};
		\node[right]  at(0,-4) (i) { $\mathcal {HQH}(L)$};
		\node[right] at(5,-4) (j) { $\mathcal {QHH}(L)$};
		\draw[->](h)--(i);
		\draw[->](i)--(j);
		\node[right] at(-2.2,-3.7)  {$\mathcal H\phi_L$};
		\node[right]  at(2.8,-3.7)  { $\phi_{\mathcal H(L)}$};
		\node[right]  at(-4.8,-6.3) (k) { $\mathcal{HQ}(L)$};
		\node[right]  at(5.2,-6.3) (l) { $\mathcal{QH}(L)$};
		\draw[->](h)--(k);
		\draw[->](j)--(l);
		\draw[->](k)--(l);
		\node[right]  at(-5.2,-5.2)  { $\mu_{\mathcal Q(L)}$};
		\node[right]  at(6.1,-5.2)  { $\mathcal Q\mu_L$};
		\node[right]  at(0.4,-6)  { $\phi_L$};
	\end{tikzpicture}
\end{figure}
\\
Let $\mathbf A$ be an arbitrary element in $\mathcal{HHQ}(L)$. By computation, we have 
\begin{center}
$\phi_L\circ \mu_{\mathcal Q(L)} (\mathbf A) = \{A\in \mathcal H(L): \forall K\in \bigcup \mathbf A, A\cap K\neq\emptyset\}$,
	
	$\mathcal Q\mu_L\circ\phi_{\mathcal H(L)}\circ\mathcal H\phi_L(\mathbf A) = {\uparrow}\{\bigcup\mathcal C: \mathcal C\in\mathcal{HH}(L), \forall \mathcal K\in cl(\{\phi_L(\mathcal A): \mathcal A\in \mathbf A\}), \mathcal C\cap \mathcal K\neq\emptyset\}$.
\end{center}
For each $A\in \phi_L\circ \mu_{\mathcal Q(L)} (\mathbf A)$, we have ${\downarrow}\{A\}\in \mathcal {HH}(L)$ and  $A = \bigcup{\downarrow}\{A\}$. Assume, for contradiction, that there exists a $\mathcal K\in cl(\{\phi_L(\mathcal A): \mathcal A\in \mathbf A\})$ such that ${\downarrow}\{A\}\cap\mathcal K = \emptyset$. Then $\mathcal K\in \Box(\mathcal H(L)\setminus{\downarrow}\{A\})$. Since $\Sigma\mathcal H(L)$ is well-filtered, $\Box(\mathcal H(L)\setminus{\downarrow}\{A\})$ is a Scott open subset of $\mathcal {QH}(L)$. Therefore, $\Box(\mathcal H(L)\setminus{\downarrow}\{A\})\cap \{\phi_L(\mathcal A): \mathcal A\in \mathbf A\}\neq\emptyset$. Hence, there exists an $\mathcal A\in \mathbf A$ such that $\phi_L(\mathcal A)\subseteq \mathcal H(L)\setminus{\downarrow}\{A\}$, that is, $\phi_L(\mathcal A)\cap {\downarrow}\{A\} = \emptyset$.
Meanwhile, this indicates $A\notin \phi_L(\mathcal A)$. It follows  that $A\cap K = \emptyset$ for some $K\in\mathcal A$.  Since $\mathcal A\in \mathbf A$ implies  $K\in \bigcup \mathbf A$, the relation $A\cap K = \emptyset$ directly contradicts  the fact that $A\in \phi_L\circ \mu_{\mathcal Q(L)} (\mathbf A)$. So $A = \bigcup{\downarrow}\{A\}\in \mathcal Q\mu_L\circ\phi_{\mathcal H(L)}\circ\mathcal H\phi_L(\mathbf A)$; consequently, 	$\phi_L\circ \mu_{\mathcal Q(L)} (\mathbf A)\subseteq \mathcal Q\mu_L\circ\phi_{\mathcal H(L)}\circ\mathcal H\phi_L(\mathbf A)$.

Pick an arbitrary element $\bigcup\mathcal C\in \mathcal Q\mu_L\circ\phi_{\mathcal H(L)}\circ\mathcal H\phi_L(\mathbf A)$, where $\mathcal C\in\mathcal{HH}(L)$ and for every $\mathcal K\in cl(\{\phi_L(\mathcal A): \mathcal A\in \mathbf A\}), \mathcal C\cap \mathcal K\neq\emptyset$.
	 Assume that there exists a $K\in \bigcup\mathbf A$, i.e., $K\in \mathcal A$ for some $\mathcal A\in \mathbf A$, such that $\bigcup\mathcal C\cap K = \emptyset$. Then for each $C\in \mathcal C$, $C\cap K = \emptyset$. This immediately implies that $C\notin\phi_L(\mathcal A)$ for all $C\in \mathcal C$; in other words, $\mathcal C\cap \phi_L(\mathcal A) = \emptyset$. This  contradicts the fact that $\bigcup\mathcal C\in \mathcal Q\mu_L\circ\phi_{\mathcal H(L)}\circ\mathcal H\phi_L(\mathbf A)$. So 	$\mathcal Q\mu_L\circ\phi_{\mathcal H(L)}\circ\mathcal H\phi_L(\mathbf A)\subseteq \phi_L\circ \mu_{\mathcal Q(L)} (\mathbf A)$.

To sum up, $\phi$ defined in Definition \ref{phi} is a distributive law of $\mathcal Q$ over $\mathcal H$.
\qedhere
\end{proof}
\end{prop}
	
\noindent
From the above results, we derive that $\mathcal Q\mathcal H$ is a  monad on $\mathbf{WF}$. 
 Notably, S. Vickers  studied the Eilenberg-Moore algebras of the double powerlocale $P_U P_L$ (where $P_U$ and $P_L$ denote the upper and lower powerlocales, which commute up to isomorphisms via a distributive law) and  concluded that $P_U P_L$-algebras correspond to localic frames \cite{steven}. Motivated by this work, we now investigate the  $\mathcal{Q}\mathcal{H}$-algebra  structures.

We first give the unit $\gamma$ and the multiplication $\rho$ of $\mathcal Q\mathcal H$ with the following diagrams, which are closely related to that of $\mathcal Q$ and $\mathcal H$ according to J. Beck \cite{beck}.
\begin{figure}[ht]
	\centering
	\begin{tikzpicture}[line width=0.6pt,scale=0.7]
		\node[right] at(-5.5,0)  { $\gamma:$};
		\node[right] at(-4,0) (a) { $Id$};
		\node[right] at(-1,0.8) (b) { $\mathcal Q$};
		\node[right]  at(-1,-0.8) (c) {$\mathcal H$};
		\draw[->](a)--(b);
		\draw[->](a)--(c);
		\node[right] at(2,0) (d) { $\mathcal Q\mathcal H$};
\draw[->](b)--(d);
\draw[->](c)--(d);
\node[right]  at(-2.5,0.7)  { $\theta$};
\node[right] at(-2.5,-0.7)  {$\eta$};
\node[right]  at(0.5,0.7)  {$\mathcal Q\eta$};
\node[right] at(0.5,-0.8)  { $\theta_{\mathcal H}$};

	\node[right]  at(-5.5,-3)  {$\rho:$};
\node[right]  at(-4,-3) (e) { $\mathcal Q\mathcal H\mathcal Q\mathcal H$};
\node[right]  at(0.2,-3) (f) {$\mathcal Q\mathcal Q\mathcal H\mathcal H$};
\node[right]  at(4.4, -2.2) (g) { $\mathcal Q\mathcal H\mathcal H$};
\node[right] at(4.4, -3.8) (h) { $\mathcal Q\mathcal Q\mathcal H$};
\node[right] at(8.4, -3) (i) {$\mathcal Q\mathcal H$};
\draw[->](e)--(f);
\draw[->](f)--(g);
\draw[->](f)--(h);
\draw[->](g)--(i);
\draw[->](h)--(i);
\node[right] at(-1.6,-2.6)  { $\mathcal Q\phi_{\mathcal H}$};
\node[right] at(2.6,-2.3)  { $\iota_{\mathcal H\mathcal H}$};
\node[right]  at(2.5,-3.8)  { $\mathcal Q\mathcal Q\mu$};
\node[right] at(6.6,-2.3)  {$\mathcal Q\mu$};
\node[right] at(6.8,-3.7)  { $\iota_{\mathcal H}$};
\end{tikzpicture}
\end{figure}

\noindent
Before  investigating the $\mathcal Q\mathcal H$-algebras, it is necessary to introduce the concept of \emph{$\rm{KZ}$-monad}, which is specialized by M. H. Escard$\acute{\rm o}$ in \cite{escardo} from A. Kock’s notion of $\rm{KZ}$-doctrine in a 2-category \cite{kock} to a poset-enriched category, where a \emph{poset-enriched category} is a category whose hom-sets are posets and the composition operator is monotone. What is more, a functor  between poset-enriched categories is called a \emph{poset-functor} if it
is monotone on hom-posets.

\begin{defiC}[\cite{escardo}]
	Let $\mathbf C$ be a poset-enriched category. A monad $\mathcal T = (\mathcal T, \eta, \mu)$ on $\mathbf C$ is called a \emph{right $\rm{KZ}$-monad}  if $\mathcal T$ is a poset-functor and satisfies $\eta_{\mathcal TX} \leq\mathcal  T\eta_{X}$ for all $X\in \mathbf C$. \emph{Left $\rm{KZ}$-monads} are defined poset-dually, by reversing the inequality.
\end{defiC}
\begin{rem} \label{kzmonad}
For a right $\rm{KZ}$-monad $\mathcal T$,	M. H. Escard$\acute{\rm o}$ in \cite{escardo} also provided some equivalent characterizations:
	\begin{enumerate}[($\mathrm{KZ}_1$)]
		\item For all $C\in \mathbf C$, an arrow $\alpha: \mathcal TC\rightarrow C$ is a structure map iff $\eta_{C}\dashv\alpha$ and $\alpha\circ\eta_{C} = id_{C}$.
		\item  $\eta_{\mathcal T C} \dashv \mu_C$ for all $C \in \mathbf C$.
		\item $\mu_{C}\dashv\mathcal T\eta_{C}$ for all $C\in \mathbf C$.
	\end{enumerate}
\end{rem}

\begin{propC}[\cite{escardo}]\label{struc}
Let $\mathcal T = (\mathcal T, \eta, \mu)$ be a $\rm{KZ}$-monad on a category $\mathbf{C}$. Then by $(\mathrm{KZ}_1)$, every object has at most one structure map.
	
\end{propC}
%

\begin{lem}
	The category $\mathbf{WF}$ is a poset-enriched category; the functors $\mathcal H$ and $\mathcal Q$ are both poset-functors.
\end{lem}

\begin{prop} \label{qhkz}
	$\mathcal H$ is a left $\rm{KZ}$-monad and $\mathcal Q$ is a right $\rm{KZ}$-monad.
	\begin{proof}
	The definitions of the two $\rm{KZ}$-monads make the proof trivial.
	\qedhere
	\end{proof}
\end{prop}

\begin{prop}\label{qhalgebra}
	
Let $(L, \alpha)$ be a $\mathcal Q\mathcal H$-algebra in $\mathbf{WF}$. Then
		\begin{enumerate}[(1)]
			
\item $(\mathcal Q\mathcal H(L), \iota_{\mathcal H(L)})$ and $(L, \alpha^\mathcal Q)$ with $\alpha^\mathcal Q = \alpha \circ \mathcal Q\eta_L$ are $\mathcal Q$-algebras, and $\alpha: (\mathcal Q\mathcal H(L), \iota_{\mathcal H(L)}) \rightarrow (L, \alpha^\mathcal Q)$ is a $\mathcal Q$-algebra homomorphism.
	
\item $(\mathcal Q\mathcal H(L), \rho_L \circ\theta_{\mathcal H(L)})$ and $(L, \alpha^\mathcal H)$ with $\alpha^\mathcal H = \alpha\circ \theta_{\mathcal H(L)}$ are $\mathcal H$-algebras and $\alpha: (\mathcal Q\mathcal H(L), \rho_L \circ \theta_{\mathcal H(L)}) \rightarrow (L, \alpha^\mathcal H)$ is an $\mathcal H$-algebra homomorphism.
	
\item $\alpha$ is the unique $\mathcal Q\mathcal H$-structure map of $L$.
		\end{enumerate}
\begin{proof}
	
	By (1) and (2) of Proposition \ref{mp}, we know $(\mathcal Q\mathcal H(L), \rho_L)$ is a $\mathcal Q\mathcal H$-algebra and $\alpha: (\mathcal Q\mathcal H(L), \rho_L) \rightarrow (L, \alpha)$ is a $\mathcal Q\mathcal H$-algebra homomorphism.
	\begin{enumerate}[(1)]

\item Theorem \ref{beck1} and Proposition \ref{dis} imply that $\mathcal Q\eta: \mathcal Q \rightarrow \mathcal Q\mathcal H$ is a morphism of monads. By Lemma \ref{algebra}, $(\mathcal Q\mathcal H(L), \rho_L \circ \mathcal Q\eta_L)$ and $(L, \alpha \circ \mathcal Q\eta_L) = (L, \alpha^\mathcal Q)$ are $\mathcal Q$-algebras, and $\alpha: (\mathcal Q\mathcal H(L), \rho_L \circ \mathcal Q\eta_L) \rightarrow (L, \alpha^\mathcal Q)$ is a $\mathcal Q$-algebra homomorphism.   By Proposition \ref{mp} (1), $(\mathcal Q\mathcal H(L), \iota_{\mathcal H(L)})$ is also a $\mathcal Q$-algebra. Since
	$\mathcal Q$ is a right $\rm{KZ}$-monad,  Proposition \ref{struc} indicates that structure maps for a $\rm{KZ}$-monad are unique, $\rho_L \circ\mathcal  Q\eta_L = \iota_{\mathcal H(L)}$ holds;  therefore, $\alpha: (\mathcal Q\mathcal H(L), \iota_{\mathcal H(L)}) \rightarrow (L, \alpha^\mathcal Q)$ is a $\mathcal Q$-algebra homomorphism.

	\item Theorem \ref{beck1} and Proposition \ref{dis} imply that $\theta_\mathcal H: \mathcal H \rightarrow \mathcal Q\mathcal H$ is a morphism of monads. By Lemma \ref{algebra}, $(\mathcal Q\mathcal H(L), \rho_L \circ \theta_{\mathcal H(L)})$ and $(L, \alpha \circ \theta_{\mathcal H(L)}) = (L, \alpha^\mathcal H)$ are $\mathcal H$-algebras, and $\alpha: (\mathcal Q\mathcal H(L), \rho_L \circ \theta_{\mathcal H(L)}) \rightarrow (L, \alpha^\mathcal Q)$ is an $\mathcal H$-algebra homomorphism.

\item   Assume that $\beta: \mathcal {QH}(L)\rightarrow L$ is also a structure map of $L$. Then $\beta^{\mathcal Q} = \beta\circ\mathcal Q\eta_{L}$ and $\beta^{\mathcal H} = \beta\circ\theta_{\mathcal H(L)}$ are $\mathcal Q$-algebra and $\mathcal H$-algebra structure maps, respectively. Since both $\mathcal Q$ and $\mathcal H$ are  $\rm{KZ}$-monads,  it follows from Proposition \ref{struc} that each $\mathcal Q$-algebra ($\mathcal H$-algebra) has at most one structure map. 	Thus $\beta^\mathcal Q = \alpha^\mathcal Q$ and $\beta^\mathcal H = \alpha^\mathcal H$. By (1), $\alpha$ and $\beta$ are $\mathcal Q$-algebra homomorphisms from $(\mathcal Q\mathcal H(L), \iota_{\mathcal H(L)})$ to $(L, \alpha^\mathcal Q) = (L, \beta^\mathcal Q)$. Since $\alpha \circ \theta_{\mathcal H(L)} = \alpha^\mathcal H = \beta^\mathcal H = \beta \circ \theta_{\mathcal H(L)}$, Proposition \ref{mp} (3) implies $\alpha = \beta$.
	
		\end{enumerate}

\end{proof}
\end{prop}	
	
\begin{prop}
	If  $L$ is a $\mathcal {QH}$-algebra, then $L$ is a frame.
	\begin{proof}
	Since $L$ is both   a $\mathcal {Q}$-algebra and  an $\mathcal {H}$-algebra by  Proposition \ref{qhalgebra}, it follows from Corollary \ref{halgebra} and Proposition \ref{qalgebra} that $L$ is  a meet continuous complete lattice.
 Thus to complete the proof, it suffices to show that $L$ is distributive.
	
	$\mathbf{Claim~1}$: $\mathcal {QH}(L)$ is a distributive lattice.
	
Obviously,	the union of arbitrary finite elements in $\mathcal {QH}(L)$ remains compact saturated. By Lemma \ref{coherent}, we know the complete lattice $\mathcal {H}(L)$ endowed with the Scott topology is coherent, that is, any finite intersection of compact saturated subsets of  $\mathcal {H}(L)$ is again a member  of $\mathcal {QH}(L)$. Consequently, $\mathcal {QH}(L)$ forms a lattice, where the supremum and infimum  of each finite subset are their  intersection and union, respectively. This implies the distributivity of $\mathcal {QH}(L)$.
	
		$\mathbf{Claim~2}$: The structure map $\alpha: \mathcal {QH}(L)\rightarrow L$ is a lattice homomorphism.

	By Proposition \ref{qhalgebra}, $\alpha$ is both a $\mathcal Q$-homomorphism and an $\mathcal H$-homomorphism. Thus $\alpha$ preserves finite infs by Lemma \ref{qhom}  and all sups by Lemma  \ref{ha}.


Now we consider the functions $\alpha\circ\varphi_{i}\circ(\theta_{\mathcal H}\circ\eta_{L})^{3}: L\times L\times L\rightarrow L,\ i = 1,2$, where $(\theta_{\mathcal H}\circ\eta_{L})^{3}: L\times L\times L\rightarrow \mathcal{QH}(L)\times \mathcal{QH}(L)\times\mathcal{QH}(L) $ is defined as
\begin{center}
	 $(\theta_{\mathcal H}\circ\eta_{L})^{3}(x, y, z) = ({\uparrow}\{{\downarrow}x\}, {\uparrow}\{{\downarrow}y\}, {\uparrow}\{{\downarrow}z\})$, for any $x, y, z\in L$,
\end{center}
and $\varphi_{i}:(\mathcal{QH}(L))^{3}\rightarrow\mathcal{QH}(L),\ i = 1,2$ are defined below:
for each triple $(\mathcal K_{1}, \mathcal K_{2}, \mathcal K_{3})\in (\mathcal{QH}(L))^{3}$,
\begin{center}
	$\varphi_{1}(\mathcal K_{1}, \mathcal K_{2}, \mathcal K_{3}) = \mathcal K_{1}\wedge (\mathcal K_{2}\vee \mathcal K_{3})$, \\
		$\varphi_{2}(\mathcal K_{1}, \mathcal K_{2}, \mathcal K_{3}) = (\mathcal K_{1}\wedge\mathcal K_{2})\vee (\mathcal K_{1}\wedge\mathcal K_{3})$.
\end{center}

\noindent
One can calculate that  $\alpha\circ\varphi_{1}\circ(\theta_{\mathcal H}\circ\eta_{L})^{3}(x, y,  z) = x\wedge(y\vee z)$  and $\alpha\circ\varphi_{2}\circ(\theta_{\mathcal H}\circ\eta_{L})^{3}(x, y, z) = (x\wedge y)\vee(x\wedge z)$.  It is obvious that $\varphi_{1} = \varphi_{2}$ by Claim 1. Thus $x\wedge(y\vee z) = (x\wedge y)\vee(x\wedge z)$, that is, $L$ satisfies the distributive law. In conclusion, $L$ is a frame.
\qedhere
	\end{proof}
	
\end{prop}	
	
\begin{prop}
	Every $\mathcal {QH}$-algebra homomorphism is a frame homomorphism.
	\begin{proof}
		Let $f: (L, \alpha)\rightarrow (M, \beta)$ be a $\mathcal {QH}$-algebra homomorphism. It is shown in Theorem \ref{beck1} that $\mathcal Q\eta: \mathcal Q\rightarrow\mathcal {QH}$ and $\theta_{\mathcal H}: \mathcal H\rightarrow\mathcal {QH}$ are morphisms of monads; hence, by Lemma \ref{algebra}, $f$ induces a $\mathcal Q$-algebra homomorphism $f: (L, \alpha\circ\mathcal Q\eta_{L})\rightarrow (M, \beta\circ\mathcal Q\eta_{M})$    and an  $\mathcal H$-algebra homomorphism $f: (L, \alpha\circ\theta_{\mathcal H(L)})\rightarrow (M, \beta\circ\theta_{\mathcal H(M)})$. Consequently, $f$ preserves  finite infs and arbitrary sups; in other words, $f$ is a frame homomorphism.
	\qedhere
	\end{proof}
\end{prop}	

\noindent
Let $\mathbf{Frm}$ be the category of all frames and all frame homomorphisms.	
	
\begin{thm}
	The Eilenberg-Moore category of $\mathcal {QH}$ over $\mathbf{WF}$ is a subcategory of the category  $\mathbf{Frm}$.
\end{thm}

\begin{rem}
	See   again the deflationary semilattice $\mathcal L$  in $\mathbf{WF}$ quoted in Example \ref{l}, it has been proven not to be a $\mathcal Q$-algebra. Thus as evidenced by Proposition \ref{qhalgebra},  $\mathcal L$ is not a $\mathcal {QH}$-algebra, which indicates that $\mathcal {QH}$ cannot give the free deflationary semilattices for $\mathbf{WF}$.
\end{rem}

\begin{rem}
The observation has been made by S. Vickers and also P. Taylor in a different context. S. Vickers has proved that the composite of the power constructions is equivalent
	to double exponentiation for exponentiable (i.e., locally compact) locales \cite{steven}, and P. Taylor characterized the algebras of the double exponentiation monad on locally compact spaces as frames \cite{taylor}.
	Furthermore, the characterizations work in more general categories when the necessary exponentials exist.
		 We know every core-compact well-filtered dcpo is locally compact sober \cite{cwl} and by Corollary \ref{lsi}, for each locally compact sober dcpo, the Hoare and the Smyth power constructions commute with each other. Thus, our composite monad $\mathcal{Q}\mathcal{H}$  is also equivalent to the monad defined as double exponentiation by the Sierpi\'nski space in the core-compact case. 
	Nevertheless, it remains a natural open question whether this equivalence can be extended to a broader category of well-filtered dcpos satisfying conditions (i) and (ii) of Theorem 4.9. We leave this for future work.
	
\end{rem}
\section*{Acknowledgment}
We would like to thank the anonymous reviewers for the professional comments and valuable suggestions.

\section*{Declaration of competing interest}
The authors declare that they have no known competing financial interests or personal relationships that could have appeared to influence the work reported in this paper.


		
		\bibliographystyle{alphaurl}
		
\bibliography{commutativity-lmcs}

\end{document}